\documentclass[a4paper]{article}

\usepackage[english]{babel}
\usepackage[utf8]{inputenc}
\usepackage{amsfonts}
\usepackage{amsmath}
\usepackage{amssymb}

\DeclareMathOperator*{\argmin}{arg\,min}

\usepackage{mathtools}
\usepackage{graphicx}
\usepackage{tikz-cd}
\usepackage{hyperref}
\hypersetup{
    colorlinks=true,
    linkcolor=red,
    filecolor=green,
    urlcolor=blue,
}

\usepackage{geometry}
\newcommand{\X}{\mathcal{X}}
\newcommand{\Y}{\mathcal{Y}}
\newcommand{\Z}{\mathcal{Z}}
\newcommand{\N}{\mathcal{N}}
\newcommand{\x}{\mathbf{x}}

\newcommand{\z}{\mathbf{z}}
\newcommand{\R}{\mathbf{R}}
\newcommand{\RR}{\mathbf{r}}
\newcommand{\T}{\mathbf{t}}
\newcommand{\n}{\mathbf{n}}
\newcommand{\I}{\mathbf{I}}
\newcommand{\0}{\mathbf{0}}
\newcommand{\A}{\mathbf{A}}
\newcommand{\bb}{\mathbf{b}}
\newcommand{\XX}{\mathbf{X}}
\newcommand{\LL}{\mathbf{L}}

\newcommand{\mycomment}[1]{}

\title {\textbf{{\huge Nonrigid Surface Registration: Gradient Calculations}}}

\author{\\[2ex] {\LARGE Dániel Unyi} \\[2ex]
Department of Telecommunications and Media Informatics \\
Budapest University of Technology and Economics \\
Műegyetem rkp. 3., H-1111 Budapest, Hungary \\[8ex]}

\date{}

\begin{document}

\maketitle

\section{Rigid registration}

The goal of rigid registration is to align a source surface $\X$ to a target surface $\Y$.
The alignment process involves iteratively transforming $\X$ closer and closer to $\Y$, such that
$\X=\Z^0 \rightarrow \Z^1 \rightarrow \Z^2 \rightarrow ... \rightarrow \Z^T=\Y$.
In other words, in each iteration we seek a rotation matrix $\tilde{\R}$ and a translation vector $\tilde{\T}$ such that $\X^{t+1} = \tilde{\R} (\R^t \X + \T^t) + \tilde{\T} = \tilde{\R} \X^t + \tilde{\T}$,
which corresponds to minimizing the error function
\[ w_1 \sum_{i=1}^N || \Pi_{\Y}(\z_i^{t+1}) - \z_i^{t+1} ||_2^2 + w_2 \sum_{i=1}^N || \tilde{\R} \x_i^t + \tilde{\T} - \z_i^{t+1} ||_2^2 \]
${\x_1, \x_2, ..., \x_n}$ and ${\z_1, \z_2, ..., \z_n}$ are sets of points sampled from the surfaces $\X$ and $\Z$,
and $\Pi_{\Y}(\z_i)$ denotes the closest point on $\Y$ to the point $\z_i$.
$w_1$ and $w_2$ are simple weighting factors.
By applying two simplifications, we can trace the minimization problem back to solving a system of linear equations.
First, there is nonlinearity in $\Pi_{\Y}(\z_i^{t+1})$, so we use the previous estimate $\Pi_{\Y}(\z_i^t)$ which is constant\footnote{For efficient closest point lookup, a k-d tree has to be built in advance.}.
Second, there is also nonlinearity in $\tilde{\R}$.
Assuming small rotations, we can linearize the rotation matrices using $\cos(\cdot) \approx 1$ and $\sin(\cdot) \approx 0$:
\[\tilde{\R} =
\begin{bmatrix} \cos\alpha & -\sin\alpha & 0 \\ \sin\alpha & \cos\alpha & 0 \\ 0 & 0 & 1 \end{bmatrix}
\begin{bmatrix} \cos\beta & 0 & \sin\beta \\ 0 & 1 & 0 \\ -\sin\beta & 0 & \cos\beta \end{bmatrix}
\begin{bmatrix} 1 & 0 & 0 \\ 0 & \cos\gamma & -\sin\gamma \\ 0 & \sin\gamma & \cos\gamma \end{bmatrix}\]
\[\tilde{\R} \approx
\begin{bmatrix} 1 & -\gamma & \beta \\ \gamma & 1 & -\alpha \\ -\beta & \alpha & 1 \end{bmatrix} =
\begin{bmatrix} 1 & 0 & 0 \\ 0 & 1 & 0 \\ 0 & 0 & 1 \end{bmatrix} + \begin{bmatrix} 0 & -\gamma & \beta \\ \gamma & 0 & -\alpha \\ -\beta & \alpha & 0 \end{bmatrix}
\]
By introducing $\tilde{\RR} = \begin{bmatrix}\alpha & \beta & \gamma\end{bmatrix}$, we can substitute $\tilde{\R}$ with $\I + \tilde{\RR} \times$.
Adding Tikhonov regularization might be beneficial in ensuring that the rotations remain small.

\subsection{Gradient calculations}

We have to minimize the following function with respect to $\tilde{\RR}$, $\tilde{\T}$, and $\z_j^{t+1}$:
\[ \argmin_{\tilde{\RR}, \tilde{\T}, \z_j^{t+1}} w_1 \sum_{i=1}^N || \Pi_{\Y}(\z_i^t) - \z_i^{t+1} ||_2^2 + w_2 \sum_{i=1}^N || \x_i^t + \tilde{\RR} \times \x_i^t + \tilde{\T} - \z_i^{t+1} ||_2^2 \]

\begin{itemize}
    \item Minimizing w.r.t. $\tilde{\RR}$:
    \[ \frac{\partial}{\partial \tilde{\RR}} \sum_{i=1}^N || \x_i^t + \tilde{\RR} \times \x_i^t + \tilde{\T} - \z_i^{t+1} ||_2^2 = 0 \]
    \[ \sum_{i=1}^N \x_i^t \times (\x_i^t + \tilde{\RR} \times \x_i^t + \tilde{\T} - \z_i^{t+1}) = \0 \]
    \[ -\sum_{i=1}^N \XX_i^t \XX_i^t \tilde{\RR} + \sum_{i=1}^N \XX_i^t \tilde{\T} - \sum_{i=1}^N \XX_i^t \z_i^{t+1} = \0 \]
    \item Minimizing w.r.t. $\tilde{\T}$:
    \[ \frac{\partial}{\partial \tilde{\T}} \sum_{i=1}^N || \x_i^t + \tilde{\RR} \times \x_i^t + \tilde{\T} - \z_i^{t+1} ||_2^2 = 0 \]
    \[ \sum_{i=1}^N \x_i^t + \tilde{\RR} \times \x_i^t + \tilde{\T} - \z_i^{t+1} = \0 \]
    \[ -\sum_{i=1}^N \XX_i^t \tilde{\RR} + N \tilde{\T} - \sum_{i=1}^N \z_i^{t+1} = -\sum_{i=1}^N \x_i^t \]
    \item Minimizing w.r.t. $\z_j^{t+1}$:
    \[ \frac{\partial}{\partial \z_j^{t+1}} ( w_1 \sum_{i=1}^N || \Pi_{\Y}(\z_i^t) - \z_i^{t+1} ||_2^2 + w_2 \sum_{i=1}^N || \x_i^t + \tilde{\RR} \times \x_i^t + \tilde{\T} - \z_i^{t+1} ||_2^2 ) = 0 \]
    \[ \frac{w_1}{w_2} (\z_j^{t+1} - \Pi_{\Y}(\z_j^t)) + \z_j^{t+1} - \x_j^t - \tilde{\RR} \times \x_j^t - \tilde{\T} = \0 \]
    \[ \XX_j^t \tilde{\RR} - \tilde{\T} + \left(1 + \frac{w_1}{w_2}\right) \z_j^{t+1} = \frac{w_1}{w_2} \Pi_{\Y}(\z_j^t) + \x_j^t \]
\end{itemize}

The result is a system of 6 + 3N linear equations with the same number of unknowns.

\[ \x = \begin{bmatrix} \tilde{\RR} \\ \tilde{\T} \\ \z_1^{t+1} \\ \z_2^{t+1} \\ \vdots \\ \z_N^{t+1} \end{bmatrix} \]

\[ \bb = \begin{bmatrix} \0 \\ -\sum_{i=1}^N \x_i^t \\ \frac{w_1}{w_2} \Pi_{\Y}(\z_1^t) + \x_1^t \\ \frac{w_1}{w_2} \Pi_{\Y}(\z_2^t) + \x_2^t \\ \vdots \\ \frac{w_1}{w_2} \Pi_{\Y}(\z_N^t) + \x_N^t \end{bmatrix} \]

\[ \A = \begin{bmatrix}
\A_{r r} & \A_{t r} & \A_{z_1 r} & \A_{z_2 r} & \dots & \A_{z_N r} \\
\A_{r t} & \A_{t t} & \A_{z_1 t} & \A_{z_2 t} & \dots & \A_{z_N t} \\
\A_{r z_1} & \A_{t z_1} & \A_{z_1 z_1} & \0 & \dots & \0 \\
\A_{r z_2} & \A_{t z_2} & \0 & \A_{z_2 z_2} & \dots & \0 \\
\vdots & \vdots & \vdots & \vdots & \ddots & \vdots \\
\A_{r z_N} & \A_{t z_N} & \0 & \0 & \dots & \A_{z_N z_N}
\end{bmatrix} \]

where
\[ \A_{r r} = -\sum_{i=1}^N \XX_i^t \XX_i^t \]
\[ \A_{t r} = \sum_{i=1}^N \XX_i^t \]
\[ \A_{z_i r} = -\begin{bmatrix} \XX_1^t & \XX_2^t & \dots & \XX_N^t \end{bmatrix} \]
\[ \A_{r t} = -\sum_{i=1}^N \XX_i^t \]
\[ \A_{t t} = N \I \]
\[ \A_{z_i t} = - \begin{bmatrix} \I & \I & \dots & \I \end{bmatrix} \]
\[ \A_{r z_j} = \begin{bmatrix} \XX_1^t \\ \XX_2^t \\ \dots \\ \XX_N^t \end{bmatrix} \]
\[ \A_{t z_j} = - \begin{bmatrix} \I \\ \I \\ \dots \\ \I \end{bmatrix} \]
\[ \A_{z_i z_j} = \left(1 + \frac{w_1}{w_2}\right) \begin{bmatrix}
\I & \0 & \dots & \0 \\
\0 & \I & \dots & \0 \\
\vdots & \vdots & \ddots & \vdots \\
\0 & \0 & \dots & \I \\                
\end{bmatrix} \]


\section{As-rigid-as-possible registration}

In each iteration we also compute the local rotation matrices $\tilde{\R}_i$:
\[ w_1 \sum_{i=1}^N || \Pi_{\Y}(\z_i^{t+1}) - \z_i^{t+1} ||_2^2 + w_2 \sum_{i=1}^N || \tilde{\R} \x_i^t + \tilde{\T} - \z_i^{t+1} ||_2^2
+ w_3 \sum_{i=1}^N \sum_{k \in \N(i)} || \tilde{\R}_i (\x_k^t - \x_i^t) - (\z_k^{t+1} - \z_i^{t+1}) ||_2^2 \]
We apply the same linearization to the local rotation matrices as we do to the global rotation matrix.
Connectivity has to be defined among the points of the source surface, resulting a graph (or a mesh specifically) with graph Laplacian $\LL$.
The neighbours of point $i$ are denoted by $\N(i)$.

\subsection{Gradient calculations}

We have to minimize the following function with respect to $\tilde{\RR}$, $\tilde{\T}$, $\tilde{\RR}_j$, and $\z_j^{t+1}$:
\begin{align*} \argmin_{\tilde{\RR}, \tilde{\T}, \tilde{\RR}_j, \z_j^{t+1}} &w_1 \sum_{i=1}^N || \Pi_{\Y}(\z_i^{t+1}) - \z_i^{t+1} ||_2^2 + w_2 \sum_{i=1}^N || \x_i^t + \tilde{\RR} \times \x_i^t + \tilde{\T} - \z_i^{t+1} ||_2^2
+ \\ + &w_3 \sum_{i=1}^N \sum_{k \in \N(i)} || (\x_k^t - \x_i^t) + \tilde{\RR}_i \times (\x_k^t - \x_i^t) - (\z_k^{t+1} - \z_i^{t+1}) ||_2^2 \end{align*}

\begin{itemize}
    \item Minimizing w.r.t. $\tilde{\RR}$:
    \[ \frac{\partial}{\partial \tilde{\RR}} \sum_{i=1}^N || \x_i^t + \tilde{\RR} \times \x_i^t + \tilde{\T} - \z_i^{t+1} ||_2^2 = 0 \]
    \[ \sum_{i=1}^N \x_i^t \times (\x_i^t + \tilde{\RR} \times \x_i^t + \tilde{\T} - \z_i^{t+1}) = \0 \]
    \[ -\sum_{i=1}^N \XX_i^t \XX_i^t \tilde{\RR} + \sum_{i=1}^N \XX_i^t \tilde{\T} - \sum_{i=1}^N \XX_i^t \z_i^{t+1} = \0 \]
    \item Minimizing w.r.t. $\tilde{\T}$:
    \[ \frac{\partial}{\partial \tilde{\T}} \sum_{i=1}^N || \x_i^t + \tilde{\RR} \times \x_i^t + \tilde{\T} - \z_i^{t+1} ||_2^2 = 0 \]
    \[ \sum_{i=1}^N \x_i^t + \tilde{\RR} \times \x_i^t + \tilde{\T} - \z_i^{t+1} = \0 \]
    \[ -\sum_{i=1}^N \XX_i^t \tilde{\RR} + N \tilde{\T} - \sum_{i=1}^N \z_i^{t+1} = -\sum_{i=1}^N \x_i^t \]
    \item Minimizing w.r.t. $\tilde{\RR}_j$:
    \[ \frac{\partial}{\partial \tilde{\RR}_j} \sum_{i=1}^N \sum_{k \in \N(i)} || (\x_k^t - \x_i^t) + \tilde{\RR}_i \times (\x_k^t - \x_i^t) - (\z_k^{t+1} - \z_i^{t+1}) ||_2^2 = 0 \]
    \[ \sum_{k \in \N(j)} (\x_k^t - \x_j^t) \times ((\x_k^t - \x_j^t) + \tilde{\RR}_j \times (\x_k^t - \x_j^t) - (\z_k^{t+1} - \z_j^{t+1})) = \0 \]
    \[ -\sum_{k \in \N(j)} (\XX_k^t - \XX_j^t) (\XX_k^t - \XX_j^t) \tilde{\RR}_j + \sum_{k \in \N(j)} (\XX_k^t - \XX_j^t) (\z_j^{t+1} - \z_k^{t+1}) = \0 \]
    \item Minimizing w.r.t. $\z_j^{t+1}$:
    \begin{align*} \frac{\partial}{\partial \z_j^{t+1}} (&w_1 \sum_{i=1}^N || \Pi_{\Y}(\z_i^t) - \z_i^{t+1} ||_2^2 + w_2 \sum_{i=1}^N || \x_i^t + \tilde{\RR} \times \x_i^t + \tilde{\T} - \z_i^{t+1} ||_2^2 + \\ + &w_3 \sum_{i=1}^N \sum_{k \in \N(i)} || (\x_k^t - \x_i^t) + \tilde{\RR}_i \times (\x_k^t - \x_i^t) - (\z_k^{t+1} - \z_i^{t+1}) ||_2^2) = 0 \end{align*}
    \begin{align*} &w_1 (\z_j^{t+1} - \Pi_{\Y}(\z_j^t)) + w_2 (\z_j^{t+1} - \x_j^t - \tilde{\RR} \times \x_j^t - \tilde{\T}) + \\ + &w_3 \sum_{k \in \N(j)} 2(\x_k^t - \x_j^t) + (\tilde{\RR}_k + \tilde{\RR}_j) \times (\x_k^t - \x_j^t) - 2(\z_k^{t+1} - \z_j^{t+1}) = \0 \end{align*}
    \begin{align*} (&w_1 + w_2) \z_j^{t+1} - 2 w_3 \sum_{k \in \N(j)} (\z_k^{t+1} - \z_j^{t+1}) + w_2 \XX_j^t \tilde{\RR} - w_2 \tilde{\T} - w_3 \sum_{k \in \N(j)} (\XX_k^t - \XX_j^t) (\tilde{\RR}_j + \tilde{\RR}_k) = \\ = &w_1 \Pi_{\Y}(\z_j^t) + w_2 \x_j^t - 2 w_3 \sum_{k \in \N(j)} (\x_k^t - \x_j^t) \end{align*}
\end{itemize}

The result is a system of 6 + 6N linear equations with the same number of unknowns.

\[ \x = \begin{bmatrix} \tilde{\RR} \\ \tilde{\T} \\ \tilde{\RR}_1^{t+1} \\ \tilde{\RR}_2^{t+1} \\ \vdots \\ \tilde{\RR}_N^{t+1} \\ \z_1^{t+1} \\ \z_2^{t+1} \\ \vdots \\ \z_N^{t+1} \end{bmatrix} \]

\[ \bb = \begin{bmatrix} \0 \\ -\sum_{i=1}^N \x_i^t \\ \0 \\ \0 \\ \vdots \\ \0 \\ w_1 \Pi_{\Y}(\z_1^t) + w_2 \x_1^t + 2 w_3 \sum_{k=1}^N \LL_{1k} \x_k^t \\ w_1 \Pi_{\Y}(\z_2^t) + w_2 \x_2^t + 2 w_3 \sum_{k=1}^N \LL_{2k} \x_k^t \\ \vdots \\ w_1 \Pi_{\Y}(\z_N^t) + w_2 \x_N^t + 2 w_3 \sum_{k=1}^N \LL_{Nk} \x_k^t \end{bmatrix} \]

\[ \A = \begin{bmatrix}
\A_{r r} & \A_{t r} & \0 & \0 & \dots & \0 & \A_{z_1 r} & \A_{z_2 r} & \dots & \A_{z_N r} \\
\A_{r t} & \A_{t t} & \0 & \0 & \dots & \0 & \A_{z_1 t} & \A_{z_2 t} & \dots & \A_{z_N t} \\
\0 & \0 & \A_{r_1 r_1} & \0 & \cdots & \0 & \A_{z_1 r_1} & \A_{z_2 r_1} & \cdots & \A_{z_N r_1} \\
\0 & \0 & \0 & \A_{r_2 r_2} & \cdots & \0 & \A_{z_1 r_2} & \A_{z_2 r_2} & \cdots & \A_{z_N r_2} \\
\vdots & \vdots & \vdots & \vdots & \ddots & \vdots & \vdots & \vdots & \ddots & \vdots \\
\0 & \0 & \0 & \0 & \cdots & \A_{r_N r_N} & \A_{z_1 r_N} & \A_{z_2 r_N} & \cdots & \A_{z_N r_N} \\
\A_{r z_1} & \A_{t z_1} & \A_{r_1 z_1} & \A_{r_2 z_1} & \dots & \A_{r_N z_1} & \A_{z_1 z_1} & \A_{z_2 z_1} & \dots & \A_{z_N z_1} \\
\A_{r z_2} & \A_{t z_2} & \A_{r_1 z_2} & \A_{r_2 z_2} & \dots & \A_{r_N z_2} & \A_{z_1 z_2} & \A_{z_2 z_2} & \dots & \A_{z_N z_2} \\
\vdots & \vdots & \vdots & \vdots & \ddots & \vdots & \vdots & \vdots & \ddots & \vdots \\
\A_{r z_N} & \A_{t z_N} & \A_{r_1 z_N} & \A_{r_2 z_N} & \dots & \A_{r_N z_N} & \A_{z_1 z_N} & \A_{z_2 z_N} & \dots & \A_{z_N z_N} \\
\end{bmatrix} \]

where
\[ \A_{r r} = -\sum_{i=1}^N \XX_i^t \XX_i^t \]
\[ \A_{t r} = \sum_{i=1}^N \XX_i^t \]
\[ \A_{z_i r} = -\begin{bmatrix} \XX_1^t & \XX_2^t & \dots & \XX_N^t \end{bmatrix} \]
\[ \A_{r t} = -\sum_{i=1}^N \XX_i^t \]
\[ \A_{t t} = N \I \]
\[ \A_{z_i t} = - \begin{bmatrix} \I & \I & \dots & \I \end{bmatrix} \]
\[ \A_{r_i r_j} = -\begin{bmatrix}
\sum_{k=1}^N \LL_{1k} (\XX_1^t - \XX_k^t)^2 & \0 & \dots & \0 \\
\0 & \sum_{k=1}^N \LL_{2k} (\XX_2^t - \XX_k^t)^2 & \dots & \0 \\
\vdots & \vdots & \ddots & \vdots \\
\0 & \0 & \dots & \sum_{k=1}^N \LL_{Nk} (\XX_N^t - \XX_k^t)^2 \\
\end{bmatrix} \]
\[ \A_{z_i r_j} = -\begin{bmatrix}
-\sum_{k=1}^N \LL_{1k} (\XX_1^t - \XX_k^t) & \LL_{12}(\XX_1 - \XX_2) & \dots & \LL_{1N}(\XX_1 - \XX_N) \\
\LL_{21}(\XX_2 - \XX_1) & -\sum_{k=1}^N \LL_{2k} (\XX_2^t - \XX_k^t) & \dots & \LL_{2N}(\XX_2 - \XX_N) \\
\vdots & \vdots & \ddots & \vdots \\
\LL_{N1}(\XX_N - \XX_1) & \LL_{N2}(\XX_N - \XX_2) & \dots & -\sum_{k=1}^N \LL_{Nk} (\XX_N^t - \XX_k^t) \\
\end{bmatrix} \]
\[ \A_{r z_j} = w_2 \begin{bmatrix} \XX_1^t \\ \XX_2^t \\ \dots \\ \XX_N^t \end{bmatrix} \]
\[ \A_{t z_j} = - w_2 \begin{bmatrix} \I \\ \I \\ \dots \\ \I \end{bmatrix} \]
\[ \A_{r_i z_j} = -w_3 \begin{bmatrix}
\sum_{k=1}^N \LL_{1k} (\XX_1^t - \XX_k^t) & \LL_{12}(\XX_1 - \XX_2) & \dots & \LL_{1N}(\XX_1 - \XX_N) \\
\LL_{21}(\XX_2 - \XX_1) & \sum_{k=1}^N \LL_{2k} (\XX_2^t - \XX_k^t) & \dots & \LL_{2N}(\XX_2 - \XX_N) \\
\vdots & \vdots & \ddots & \vdots \\
\LL_{N1}(\XX_N - \XX_1) & \LL_{N2}(\XX_N - \XX_2) & \dots & \sum_{k=1}^N \LL_{Nk} (\XX_N^t - \XX_k^t) \\
\end{bmatrix} \]
\[ \A_{z_i z_j} = (w_1 + w_2) \begin{bmatrix}
\I & \0 & \dots & \0 \\
\0 & \I & \dots & \0 \\
\vdots & \vdots & \ddots & \vdots \\
\0 & \0 & \dots & \I \\
\end{bmatrix} + 2 w_3 \begin{bmatrix}
\LL_{11} \I & \LL_{12} \I & \dots & \LL_{1N} \I \\
\LL_{21} \I & \LL_{22} \I & \dots & \LL_{2N} \I \\
\vdots & \vdots & \ddots & \vdots \\
\LL_{N1} \I & \LL_{N2} \I & \dots & \LL_{NN} \I \\
\end{bmatrix} \]

\section{Accelerating the convergence: point-to-plane error}

In this case, we also measure the distance from the tangent plane belonging to the closest point, adding an extra term with weighting factor $w_4$:
\[ w_4 ( \n_i^T(\Pi_{\Y}(\z_i^t) - \z_i^{t+1}) )^2 \]
where $\n_i$ is the surface normal at point $i$. The gradient of this term w.r.t $\z_j^{t+1}$ is
\[ 2 w_4 ( \n_i \circ \n_i ) ( \z_i^{t+1} - \Pi_{\Y}(\z_i^t) )^2 \]



\section*{References}
\begin{enumerate}
    \item Bouaziz, S., Tagliasacchi, A., Li, H., \& Pauly, M. (2016). Modern techniques and applications for real-time non-rigid registration. In SIGGRAPH ASIA 2016 Courses (pp. 1-25).
    \item Sorkine, O., \& Alexa, M. (2007, July). As-rigid-as-possible surface modeling. In Symposium on Geometry processing (Vol. 4, pp. 109-116).
    \item Low, K. L. (2004). Linear least-squares optimization for point-to-plane icp surface registration. Chapel Hill, University of North Carolina, 4(10), 1-3.
\end{enumerate}

\end{document}